\documentstyle[12pt]{article}
\begin{document}
\centerline{\bf{The Quadratic Formula Made Hard}}
\centerline{or}
\centerline{\bf{A Less Radical Approach to Solving Equations}}
\vskip .2in
\centerline{M.L.Glasser}
\vskip .1in
\centerline{Physics Department, Clarkson University}
\centerline{Potsdam, N.Y. 13699-5280}
\vskip .3in
\centerline{\bf{Introduction}}
\vskip .1in

It appears that, along with many of my friends and colleagues, I had
been brainwashed by the great and tragic lives of Abel and Galois to
believe that no general formulas are possible for roots of equations
higher than quartic. This seemed to be confirmed by the brilliant and
arduous solution of the general quintic by Hermite. Yet, below we find
a formula giving a root to any algebraic equation of degree 2-5 and any
reduced equation (see below) of higher degree. This algorithm,
which must have been familiar to Lagrange, resulted when I was
working on a paper on the asymptotics of hypergeometric functions
where Gauss' multiplication formula for the gamma function is used
to
reduce certain infinite series, and by a happy accident my copy of
Whittaker and Watson opened at p. 133.
\vskip .1in
\centerline{\bf{ The Formula}}

   Without loss of generality it is sufficient to find at least
one root to the reduced equation
$$x^N-x+t=0 \;\;\;\;\;(N=2,3,4\dots).\eqno(1)$$
Letting  $x=\zeta^{-1/(N-1)}$,  we easily find that (1) becomes
$$\zeta=e^{2\pi i}+t\phi(\zeta)\eqno(2)$$
where
$$\phi(\zeta)=\zeta^{N/(N-1)}.\eqno(3)$$
Lagrange's theorem states that for any function $f$ analytic
in a neighborhood of a root of (2)
$$f(\zeta)=f(e^{2\pi i})+\sum_{n=1}^{\infty}\frac{t^n}{n!}
\frac{d^{n-1}}{da^{n-1}}[f^{\prime}(a)|\phi(a)|^n]_{a=e^{2\pi
i}}.\eqno(4)$$
We now  simply let $f(\zeta)=\zeta^{-1/(N-1)}$, carry out the
elementary differentiations (noting that
$D_kx^p=\Gamma(p+1)x^{p-k}/\Gamma(p-k+1)$) and we come up with
the root
$$x_1=exp[-2\pi i/(N-1)]-\frac{t}{N-1}\sum_{n=0}^{\infty}
\frac{(te^{2\pi i/(N-1)})^n}{\Gamma(n+2)}\frac{\Gamma(\frac{Nn}{N-1}+1)}{\Gamma(
\frac{n}{N-1}+1)}.\eqno(5)$$
(N-2 further roots are found by replacing $exp(2\pi i/(N-1))$ by the other
N-1-st  roots of unity, and the remaining root from the relation
$\sum x_j=\delta_{N,2}$).
By the use of Gauss' multiplication theorem, the infinite series
can be broken up into a (finite) sum of hypergeometric functions.
$$x_1=\omega^{-1}-\frac{t}{(N-1)^2}\sqrt{\frac{N}{2\pi(N-1)}}
\sum_{q=0}^{N-2}(\frac{\omega t}{N-1})^qn^{qN/(N-1)}
\frac{\prod_{k=0}^{N-1}\Gamma(\frac{Nq/(N-1)+1+k}{N})}{\Gamma
(\frac{q}{N-1}+1)\prod_{k=0}^{N-2}\Gamma(\frac{q+k+2}{N-1})}\times$$
$$\;_{N+1}F_N[\frac{qN/(N-1)+1}{N},\dots , \frac{qN/(N-1)+N}{N},1;$$
$$
\frac{q+2}{N-1},\dots ,\frac{q+N}{N-1},\frac{q}{N-1}+1;(\frac{t
\omega}{N-1})^{N-1}N^N],$$
where $\omega=exp(\frac{2\pi i}{N-1})$. In practice, $\;_{N+1}F_N$
will always be reducible to at least $\;_NF_{N-1}$. Hence the root
is a sum of at most $N-1$ hypergeometric functions.
The one technical point is that the convergence of these series
requires that t be "sufficiently small", but this can be overcome
by certain hypergeometric identities tantamount to analytic
continuation.

\vskip .2in
\centerline{\bf{Examples}}
\vskip .1in

\noindent
{\bf{N=2}}

$$x^2-x+t=0$$

Here we have 
$$x_1=1-t\sum_{n=0}^{\infty}\frac{t^n}{\Gamma(n+1)}\frac{\Gamma(2n+1)
}{\Gamma(n+1)}.\eqno(6)$$
However, by Gauss' formula
$$\Gamma(2n+1)=4^n(1/2)_n(1)_n \;\;((n)_k=\Gamma(n+k)/\Gamma(n))
\eqno(7)$$
so $$x_1=1-t\;_2F_1(1/2,1;2;4t) \eqno(8)$$
Since 
$$\;_2F_1(1/2,1;2;z)=\frac{2}{z}\{\begin{array}{cc}
1-\sqrt{1-z}&|z|\le1\\
1-i\sqrt{z-1}&|z|>1
\end{array}\eqno(8)$$
we reproduce the quadratic formula. Note that the second root
comes from
$x_1+x_2=1$.
\vskip .2in

\noindent
{\bf{N=3}}
$$x^3-x+t=0$$

By separating the sum in (5) into sums over the even and odd values of n we
obtain
$$x_1=-1+\frac{t}{2}\sum_{n=0}^{\infty}\frac{\Gamma(3n+1)t^{2n}}{
\Gamma(n+1)\Gamma(2n+2)} +\frac{t^2}{2}\sum_{n=0}^{\infty}\frac{
\Gamma(3n+5/2)t^{2n}}{\Gamma(n+3/2)\Gamma(2n+3)}.\eqno(9)$$
By breaking up the gamma functions of multiple argument by using Gauss'
multiplication theorem, the sums are easily identified as hypergeometric
series:
$$x_1=-1-\frac{t}{2}\;_2F_1(1/3,2/3;3/2;27t^2/4)+\frac{3t^2}{8}\;_3F_2(
5/6,7/6,1;3/2,2;27t^2/4).\eqno(10)$$
However, from A.P. Prudnikov et.al, Integrals and Series, Vol.3[Gordon and
Breach, 1990]) we find
$$\;_2F_1(1/3,2/3;3/2;z)=\frac{3}{\sqrt{z}}\sin(\frac{1}{3}\sin^{-1}
\sqrt{z})$$
$$\;_3F_2(5/6,7/6,1;3/2,2;z)=\frac{18}{z}[\cos(\frac{1}{3}\sin^{-1}
\sqrt{z})-1],\eqno(11)$$
and we therefore have the three roots
$$x_1=-\frac{1}{\sqrt{3}}\sin[\frac{1}{3}\sin^{-1}(t\sqrt{27}/2)]
-cos[\frac{1}{3}\sin^{-1}(t\sqrt{27}/2)]$$
$$x_2=-\frac{1}{\sqrt{3}}\sin[\frac{1}{3}\sin^{-1}(t\sqrt{27}/2)]
+\cos[\frac{1}{3}\sin^{-1}(t\sqrt{27}/2)]$$
$$x_3=\frac{2}{\sqrt{3}}\sin[\frac{1}{3}\sin^{-1}(t\sqrt{27}/2)].\eqno(12)$$
Once again, for $t>2/\sqrt{27}$ equation (10) must be analytically continued to obtain the correct form of (12). This
amounts to writing $\sin^{-1}z=\frac{\pi}{2}-iLn(z+\sqrt{z^2-1})$.
\vskip .2in
\centerline{\bf{Conclusion}}

\vskip .1in
For N=2,3,4 Eq.(5) is definitely not preferable to the standard formulas, but for N=5, e.g. we get the root
$$x=t\;_4F_3[{1/5,2/5,3/5,4/5;\atop 1/2,3/4,5/4;}\frac{3125t^4}{256}]\eqno(13)$$
in an elementary fashion with considerably less difficulty than by following the procedure in Davis' book [Introduction to
Nonlinear Ordinary Differential Equations (Dover)]. It might also
be pointed out that the above procedure carries over in a trivial way
to the trinomial equation $$y^N-ay^{N-1}+a=0\eqno(14)$$
where $y=1/x,\; a=1/t$. Numerically, these formulas are not of much use since solutions can be obtained
with the push of a button on many pocket calculators, but formulas such as (13) should have numerous entertaining uses, such as
summing the odd hypergeometric series. 

\vskip .2in
\centerline{\bf{Acknowledgements}}
\vskip .1in
This work was carried out at Melbourne University to which I am grateful for having me visit. In particular, I thank Dr. N.E. Frankel
for his money and Prof. A.J. Guttmann for his office ( with  its copy of Whittaker and Watson).

\end{document}